\documentclass[12pt,twosided,reqno]{amsart}
\usepackage{amsmath}
\usepackage{amsthm}
\usepackage{amsfonts}
\usepackage{amssymb}

\theoremstyle{theorem}

\theoremstyle{definition}

\begin{document}

\vskip .5cm
\title{Area of  triangles associated with a curve II}
 \vskip0.5cm
\thanks{
    \newline\indent { $^1$   supported by Basic Science Research Program through
    the National Research Foundation of Korea (NRF) funded by the Ministry of Education, Science and Technology (2010-0022926).}
\newline\indent { $^2$   supported by Basic Science Research Program through the National Research Foundation of Korea (NRF)
funded by the Ministry of Education, Science and Technology (2012R1A1A2042298).}}

\vskip 0.5cm

\maketitle

\vskip 0.5cm
\centerline{\scshape
Dong-Soo Kim$^1$, Wonyong Kim, Young Ho Kim$^2$ and Dae Heui Park} \vskip .2in

\begin{abstract}
It is well known that the area $U$ of the triangle
formed by three tangents to a parabola $X$ is half of the area $T$
 of the
triangle formed by joining their points of contact.
In this article, we
consider whether this property and similar ones characterizes parabolas.
As a result, we present three conditions  which are necessary and sufficient for
a strictly convex curve in the plane to be an open part of a parabola.

\end{abstract}

\vskip 1cm

\date{}
\maketitle

\section{Introduction}
 \vskip 0.50cm
Suppose that $X$ is a   regular  curve in the  plane
 ${\mathbb R}^{2}$ with nonvanishing curvature,
and   $P=A$, $A_i, i = 1, 2,$  are three distinct neighboring points on the curve $X$.
Let us denote by $\ell$, $\ell_1$, $\ell_2$ the tangent
lines passing through the points $A,A_1,A_2$ and by $B,B_1,B_2$
the intersection points $\ell_1\cap \ell_2$, $\ell\cap \ell_1$,
$\ell\cap \ell_2$, respectively.
If $X$ is an open part of a parabola, it is well known that the area $U=|\bigtriangleup BB_1B_2|$ of the triangle
formed by three tangents to the parabola $X$ is half of the area $T=|\bigtriangleup AA_1A_2|$ of the
triangle formed by joining their points of contact (\cite{D}).

A regular plane curve $X:I \rightarrow {\mathbb R}^{2}$  defined on an open interval
 is called {\it convex} if, for all $t\in I$, the trace $X(I)$ lies
 entirely on one side of the closed half-plane determined by the tangent line at $X(t)$ (\cite{dC}).

From now on,  we will  say  that a simple convex curve $X$ in the  plane
 ${\mathbb R}^{2}$ is  {\it strictly convex} if the curve   is smooth (that is, of class $C^{(3)}$) and is of positive  curvature $\kappa$
 with respect to the unit normal $N$ pointing to the convex side.
 Hence, in this case we have $\kappa(s)=\left< X''(s), N(X(s))\right> >0$, where $X(s)$ is an arclength parametrization of $X$.

For a smooth function $f:I\rightarrow {\mathbb R}$ defined on an open interval, we will also say that $f$ is
{\it strictly convex} if the graph of $f$ has  positive curvature $\kappa$ with respect to the upward unit normal $N$. This condition is equivalent to the positivity of $f''(x)$ on $I$.

Suppose that $X$ is a strictly convex  curve  in the plane
 ${\mathbb R}^{2}$ with the unit normal $N$ pointing to the convex side.
 For a fixed point $P=A \in X$, and for a sufficiently small $h>0$, we consider the  line $m$ passing through
 $P+hN(P)$ which is parallel to
 the tangent $\ell$ of $X$ at $P$ and the points  $A_1$ and $A_2$  where the line $m$ intersects the curve $X$.

We  denote by  $\ell_1$, $\ell_2$ the tangent
lines of $X$ at the points $A_1,A_2$ and by $B,B_1,B_2$
the intersection points $\ell_1\cap \ell_2$, $\ell\cap \ell_1$,
$\ell\cap \ell_2$, respectively.
We let $L_P(h)$ and  $\ell_P(h)$ denote the length $|A_1A_2|$ and $|B_1B_2|$  of the corresponding segments, respectively.

Now, we consider  $T_P(h), U_P(h)$,  $V_P(h)$ and $W_P(h)$
defined by the area  $|\bigtriangleup AA_1A_2|$,
$|\bigtriangleup BB_1B_2|$,  $|\bigtriangleup BA_1A_2|$
 of corresponding triangles and the area $|\square A_1A_2B_2B_1|$ of trapezoid $\square A_1A_2B_2B_1$, respectively.
 Then, obviously  we have
 $$T_P(h)=\frac{h}{2}L_P(h)$$ and
 $$V_P(h)=W_P(h)+U_P(h).$$

 Let us denote by $S_P(h)$  the area of the region bounded by the curve $X$ and chord $A_1A_2$.
 Then, we have  (\cite{KK4})
 $$S_P'(h)=L_P(h).$$

It is well known that  parabolas satisfy the following properties (\cite{D,S}).

\vskip 0.3cm

 \noindent {\bf Lemma 1.1.} Suppose that $X$ is an open part of a parabola.
 For arbitrary point $P\in X$ and sufficiently small $h>0$, it satisfies
  \begin{equation}\tag{1.1}
   \begin{aligned}
  S_P(h)=\frac{4}{3}T_P(h),
    \end{aligned}
   \end{equation}
 \begin{equation}\tag{1.2}
   \begin{aligned}
 S_P(h)=\frac{2}{3}V_P(h),
    \end{aligned}
   \end{equation}
   \begin{equation}\tag{1.3}
   \begin{aligned}
S_P(h)=\frac{8}{9}W_P(h),
    \end{aligned}
   \end{equation}
and
    \begin{equation}\tag{1.4}
   \begin{aligned}
U_P(h)=\frac{1}{2}T_P(h).
    \end{aligned}
   \end{equation}
 \vskip 0.3cm

 Actually, Archimedes showed that parabolas satisfy (1.1) (\cite{S}).
 Recently, in \cite{KK4} the first and third authors of the present paper proved that (1.1) is a characteristic property of parabolas
 and established some  characterizations of parabolas, which are the converses of
well-known properties of parabolas  originally due to Archimedes
(\cite{S}).  For the higher dimensional analogues of some results in \cite{KK4}, see \cite{KK2} and \cite{KK3}.

\vskip 0.3cm

 In this article,  we
study whether the remaining properties in Lemma 1.1  characterize parabolas.

First of all, in Section 3 we prove the following:
 \vskip 0.3cm

 \noindent {\bf Theorem 1.2.}
 Let $X$ denote a strictly convex $C^{(3)}$ curve  in the plane
 ${\mathbb R}^{2}$. Then the following are equivalent.

\noindent 1) For all $P\in X$ and sufficiently small $h>0$, $S_P(h)=\lambda(P)V_P(h)$,
where  $\lambda(P)$ is a function of $P$.

\noindent 2) For all $P\in X$ and sufficiently small $h>0$, $S_P(h)=\frac{2}{3}V_P(h)$.

\noindent 3) $X$ is an open part of a parabola.

\vskip 0.3cm
Next, in Section 4 we study plane curves satisfying (1.4) in Lemma 1.1.

In \cite{Kr}, Krawczyk showed that for a strictly convex $C^{(4)}$ curve $X$ in the plane
 ${\mathbb R}^{2}$, the following holds:
  \begin{equation}\tag{1.5}
   \begin{aligned}
 \lim_{A_1,A_2 \rightarrow A} \frac{T}{U}=2.
         \end{aligned}
   \end{equation}
His application of (1.5) states  that if a strictly convex $C^{(4)}$ curve $X$ in the plane
 ${\mathbb R}^{2}$ satisfies for some function $\lambda(P)$
  \begin{equation}\tag{1.6}
   \begin{aligned}
U=\lambda(P)T,
         \end{aligned}
   \end{equation}
   then $\lambda(P)=\frac{1}{2}$ and $X$ is an open part of
 the graph of a quadratic polynomial.

 Extending the results in \cite{Kr}, in \cite{KS} the first author and K.-C. Shim showed that
 if a strictly convex $C^{(3)}$ curve $X$ in the plane
 ${\mathbb R}^{2}$ satisfies (1.6) for some function $\lambda(P)$,
    then $\lambda(P)=\frac{1}{2}$ and $X$ is an open part of
 a parabola. For a further study, they also posed a question (Question 3 in \cite{KS}) whether  the property (1.4) is a characteristic property of parabolas.

\vskip 0.30cm
As a result, in Section 4  we give an affirmative answer to Question 3 in \cite{KS} as follows:

 \vskip 0.3cm
\noindent {\bf Theorem 1.3.}
 Let $X$ denote a strictly convex $C^{(3)}$ curve  in the plane
 ${\mathbb R}^{2}$. Then the following are equivalent.

\noindent 1) For all $P\in X$ and sufficiently small $h>0$, $U_P(h)=\lambda(P)T_P(h)$,
where  $\lambda(P)$ is a function of $P$.

\noindent 2) For all $P\in X$ and sufficiently small $h>0$, $U_P(h)=\frac{1}{2}T_P(h)$.

\noindent 3) $X$ is an open part of a parabola.

\vskip 0.3cm

Finally, in Section 5 we prove the following:
\vskip 0.3cm
 \noindent {\bf Theorem 1.4.}
 Let $X$ denote a strictly convex $C^{(3)}$ curve  in the plane
 ${\mathbb R}^{2}$. Then the following are equivalent.

\noindent 1) For all $P\in X$ and sufficiently small $h>0$, $S_P(h)=\lambda(P)W_P(h)$,
where  $\lambda(P)$ is a function of $P$.

\noindent 2) For all $P\in X$ and sufficiently small $h>0$, $S_P(h)=\frac{8}{9}W_P(h)$.

\noindent 3) $X$ is an open part of a parabola.
 \vskip 0.3cm

 For some characterizations of parabolas or  conic sections by  properties of tangent lines, see
\cite{KKa} and \cite{KKP}. In \cite{KK1}, using curvature function $\kappa$ and support function $h$
of a plane curve,
the first and third authors of the present paper gave a
characterization of ellipses and hyperbolas centered at the origin.

 In \cite{R}, B.  Richmond and T. Richmond  established a dozen necessary and sufficient conditions
 for the graph of a function to
be a parabola  by using elementary techniques.
\vskip 0.3cm
  Throughout this article, all curves are of class $C^{(3)}$ and connected, unless otherwise mentioned.
  \vskip 0.50cm

\section{Preliminaries }
 \vskip 0.50cm

 In order to prove Theorems in Section 1,  we  need the following lemma.
\vskip 0.3cm

 \noindent {\bf Lemma 2.1.} Suppose that   $X$  is a strictly convex  $C^{(3)}$ curve  in the plane
 ${\mathbb R}^{2}$ with  the unit normal $N$ pointing to the convex side.  Then
 we have
  \begin{equation}\tag{2.1}
   \begin{aligned}
   \lim_{h\rightarrow 0} \frac{1}{\sqrt{h}}L_P(h)= \frac{2\sqrt{2}}{\sqrt{\kappa(P)}},
    \end{aligned}
   \end{equation}
      \begin{equation}\tag{2.2}
   \begin{aligned}
   \lim_{h\rightarrow 0} \frac{1}{h\sqrt{h}}S_P(h)= \frac{4\sqrt{2}}{3\sqrt{\kappa(P)}},
    \end{aligned}
   \end{equation}
\begin{equation}\tag{2.3}
   \begin{aligned}
 \lim_{h\rightarrow 0} \frac{T_P(h)}{h\sqrt{h}}=\frac{\sqrt{2}}{\sqrt{\kappa(P)}}
         \end{aligned}
   \end{equation}
 and
 \begin{equation}\tag{2.4}
   \begin{aligned}
 \lim_{h\rightarrow 0} \frac{U_P(h)}{h\sqrt{h}}=\frac{\sqrt{2}}{2\sqrt{\kappa(P)}},
         \end{aligned}
   \end{equation}
 where $\kappa(P)$ is the curvature of $X$ at $P$ with respect to  the unit normal $N$.
 \vskip 0.3cm
 \noindent {\bf Proof.} It follows from \cite{KK4}  that (2.1) and (2.2) hold. For a proof of
 (2.3) and (2.4), see \cite{KS}. $\square$

\vskip 0.3cm
First, we obtain

 \vskip 0.3cm
\noindent {\bf Lemma 2.2.} Suppose that   $X$  is a strictly convex  $C^{(3)}$ curve  in the plane
 ${\mathbb R}^{2}$ with  the unit normal $N$ pointing to the convex side.  Then
 we have

  \begin{equation}\tag{2.5}
   \begin{aligned}
   \lim_{h\rightarrow 0} \frac{1}{\sqrt{h}}\ell_P(h)= \frac{\sqrt{2}}{\sqrt{\kappa(P)}},
    \end{aligned}
   \end{equation}

  \begin{equation}\tag{2.6}
   \begin{aligned}
   \lim_{h\rightarrow 0} \frac{1}{h\sqrt{h}}V_P(h)= \frac{2\sqrt{2}}{\sqrt{\kappa(P)}}
    \end{aligned}
   \end{equation}
   and
      \begin{equation}\tag{2.7}
   \begin{aligned}
   \lim_{h\rightarrow 0} \frac{1}{h\sqrt{h}}W_P(h)= \frac{3\sqrt{2}}{2\sqrt{\kappa(P)}}.
    \end{aligned}
   \end{equation}

   \vskip 0.3cm
 \noindent {\bf Proof.}
We  fix an arbitrary  point $P$ on $X$.
Then, we may take a coordinate system $(x,y)$
 of  ${\mathbb R}^{2}$: $P$ is taken to be the origin $(0,0)$ and $x$-axis is the tangent line $\ell$ of $X$ at $P$.
 Furthermore, we may regard $X$ to be locally  the graph of a non-negative strictly convex  function $f: {\mathbb R}\rightarrow {\mathbb R}$
 with $f(0)=f'(0)=0$. Then $N$ is the upward unit normal.

 Since the curve $X$ is of class $C^{(3)}$, the Taylor's formula of $f(x)$ is given by
 \begin{equation}\tag{2.8}
 f(x)= ax^2 + f_3(x),
  \end{equation}
where  $2a=f''(0)$ and $f_3(x)$ is an $O(|x|^3)$  function.
From  $\kappa(P)=f''(0)>0$, we see that $a$ is positive.

For a sufficiently small $h>0$, we denote by $A_1(s,f(s))$ and $A_2(t,f(t))$
 the points where the line $m:y=h$ meets the curve $X$ with $s<0<t$.
 Then $f(s)=f(t)=h$ and we get $B_1(s-h/f'(s),0)$, $B_2(t-h/f'(t),0)$ and $B(x_0,y_0)$ with
  \begin{equation}\tag{2.9}
x_0=\frac{tf'(t)-sf'(s)}{ f'(t)-f'(s)}
  \end{equation}
and
  \begin{equation}\tag{2.10}
y_0=\frac{(t-s)f'(t)f'(s)+h(f'(t)-f'(s))}{ f'(t)-f'(s)}.
  \end{equation}
Noting that $L_P(h)=t-s$ and
\begin{equation}\tag{2.11}
\ell_P(h)=t-s+\frac{f'(t)-f'(s)}{ f'(s)f'(t)}h,
  \end{equation}
 one obtains
 \begin{equation}\tag{2.12}
  \begin{aligned}
2V_P(h)&=\{t-s \}(h-y_0)\\
&=\frac{-f'(s)f'(t)}{ f'(t)-f'(s)}L_P(h)^2
   \end{aligned}
  \end{equation}
and
 \begin{equation}\tag{2.13}
  \begin{aligned}
2W_P(h)&=\{2(t-s)+ \frac{f'(t)-f'(s)}{f'(s)f'(t)}h\}h\\
&=2hL_P(h)+ \frac{f'(t)-f'(s)}{ f'(s)f'(t)}h^2.\\
   \end{aligned}
  \end{equation}

If we let
\begin{equation}\tag{2.14}
  \begin{aligned}
\alpha_P(h)= \frac{f'(t)-f'(s)}{-f'(s)f'(t)}\sqrt{h},
   \end{aligned}
  \end{equation}
then from (2.11)-(2.13) we have
\begin{equation}\tag{2.15}
  \begin{aligned}
  \ell_P(h)&=L_P(h)-\alpha_P(h)\sqrt{h},\\
2V_P(h)&=\frac{\sqrt{h}}{ \alpha_P(h)}L_P(h)^2
   \end{aligned}
  \end{equation}
and
 \begin{equation}\tag{2.16}
  \begin{aligned}
2W_P(h)&=2hL_P(h)-\alpha_P(h)h\sqrt{h}.\\
   \end{aligned}
  \end{equation}

On the other hand, from Lemma 5 in \cite{KS} we get
 \begin{equation}\tag{2.17}
  \begin{aligned}
\lim_{h\rightarrow0}\alpha_P(h)=\frac{\sqrt{2}}{\sqrt{\kappa(P)}}.
   \end{aligned}
  \end{equation}
Hence, together with Lemma 2.1, (2.15) and (2.16), this completes the proof. $\square$

\vskip 0.3cm
Next, we get the following lemma which plays a crucial role in the proofs of Theorems 1.3 and 1.4.

 \vskip 0.3cm
\noindent {\bf Lemma 2.3.} Suppose that   $X$  is a strictly convex  $C^{(3)}$ curve  in the plane
 ${\mathbb R}^{2}$ with  the unit normal $N$ pointing to the convex side.  Then
 we have
 \begin{equation}\tag{2.18}
  \begin{aligned}
\ell_P(h)=L_P(h)-h\frac{d}{dh}L_P(h).
   \end{aligned}
  \end{equation}

   \vskip 0.3cm
 \noindent {\bf Proof.}
 As in the proof of Lemma 2.2, for an arbitrary  point $P$ on $X$ we take a coordinate system $(x,y)$
 of  ${\mathbb R}^{2}$ so that (2.8) holds.
 If we put $f(t)=h$ for sufficiently small $t>0$, then the line $m:y=h$ meets the curve $X$
 at the points  $A_1(s(t),h)$ and $A_2(t,h)$ with $s=s(t)<0<t$.
 Hence we have
  \begin{equation}\tag{2.19}
  \begin{aligned}
f(s(t))=f(t)
   \end{aligned}
  \end{equation}
and
 \begin{equation}\tag{2.20}
  \begin{aligned}
B_1=(s(t)-\frac{h}{f'(s(t))},0), \quad B_2=(t-\frac{h}{f'(t)},0).
   \end{aligned}
  \end{equation}

  Thus, we get
  \begin{equation}\tag{2.21}
  \begin{aligned}
L_P(h)=t-s(t)
   \end{aligned}
  \end{equation}
and
 \begin{equation}\tag{2.22}
  \begin{aligned}
\ell_P(h)=t-s(t)-\{\frac{1}{f'(t)}-\frac{1}{f'(s(t))}\}h.
   \end{aligned}
  \end{equation}
Noting $h=f(t)$, one obtains from (2.21) that
 \begin{equation}\tag{2.23}
  \begin{aligned}
\frac{d}{dh}L_P(h)=\frac{1-s'(t)}{f'(t)}.
   \end{aligned}
  \end{equation}

Therefore, it follows from (2.19) that
\begin{equation}\tag{2.24}
  \begin{aligned}
\frac{d}{dh}L_P(h)=\frac{1}{f'(t)}-\frac{1}{f'(s(t))}.
   \end{aligned}
  \end{equation}
Together with (2.21) and (2.22), this  completes the proof. $\square$

\vskip 0.50cm

 \section{Proof of Theorem 1.2}
 \vskip 0.50cm
In this section, we prove Theorem 1.2.

 It is trivial to show that any open part of parabolas satisfy 1) and 2) in Theorem 1.2.

Conversely, suppose that  $X$ denotes a strictly convex $C^{(3)}$ curve  in the plane
 ${\mathbb R}^{2}$ which satisfies $S_P(h)=\lambda(P)V_P(h)$ for all $P\in X$ and sufficiently small $h>0$.
 Then, it follows from  Lemmas 2.1 and 2.2 that $\lambda(P)=\frac{2}{3}$.

 We  fix an arbitrary  point $A_1$ on $X$.
Then, we may take a coordinate system $(x,y)$
 of  ${\mathbb R}^{2}$ so that  $A_1$ is the origin $(0,0)$
 and $x$-axis is the tangent line of $X$ at $A_1$.
 Furthermore, we may regard $X$ to be locally
 the graph of a non-negative strictly convex  function $f: {\mathbb R}\rightarrow {\mathbb R}$
 with $f(0)=f'(0)=0$. Then $N$ is the upward unit normal.

 Since the curve $X$ is of class $C^{(3)}$, the Taylor's formula of $f(x)$ is given by
 \begin{equation}\tag{3.1}
 f(x)= ax^2 + f_3(x),
  \end{equation}
where  $2a=f''(0)$ and $f_3(x)$ is an $O(|x|^3)$  function.
From  $\kappa(A_1)=f''(0)>0$, we see that $a$ is positive.

For any point $A_2(t,f(t))$ with sufficiently small $t$,
we denote by $P=A$ the point on $X$ such that
the chord $A_1A_2$ is parallel to the tangent of $X$ at $P=A$.
Then  we have
$A=(g(t), f(g(t)))$, for a function $g: {\mathbb R}\setminus \{0\}
 \rightarrow {\mathbb R}$ which satisfies $|g(t)|<|t|$ and

\begin{equation}\tag{3.2}
   \begin{aligned}
  tf'(g(t))=f(t).
 \end{aligned}
   \end{equation}
Since   $g(t)$ tends to $0$ as $t\rightarrow 0$, we may assume that $g(0)=0$.
 \vskip 0.3cm

Then we have
\begin{equation}\tag{3.3}
   \begin{aligned}
B_1=(g(t)-\frac{tf(g(t))}{f(t)},0),  B=(t-\frac{f(t)}{f'(t)}, 0)
 \end{aligned}
   \end{equation}
   and
\begin{equation}\tag{3.4}
   \begin{aligned}
B_2=(t+\frac{tf(g(t))-f(t)g(t)}{tf'(t)-f(t)},f(t)+\frac{tf(g(t))-f(t)g(t)}{tf'(t)-f(t)}f'(t)).
 \end{aligned}
   \end{equation}

If we let $h$ the distance from $P=A$ to the chord $A_1A_2$, then by the definition of $V_P(h)$
we have
\begin{equation}\tag{3.5}
   \begin{aligned}
 V_P(h)=\epsilon \frac{f(t)}{2}\{t-\frac{f(t)}{f'(t)}\},
 \end{aligned}
   \end{equation}
   where $\epsilon =1$ if $t>0$ and $\epsilon =-1$ otherwise.
\vskip 0.3cm
We now prove the following lemma, which is useful in the proof of Theorem 1.2.

\vskip 0.3cm

 \noindent {\bf Lemma 3.1.} Suppose that $S_P(h)=\frac{2}{3}V_P(h)$. Then the function $f(t)$ satisfies
\begin{equation}\tag{3.6}
  \begin{aligned}
2f(t)^2f''(t)=f'(t)^2\{tf'(t)-f(t)\}.
   \end{aligned}
  \end{equation}
\vskip 0.3cm
 \noindent {\bf Proof.} Note that
 \begin{equation}\tag{3.7}
  \begin{aligned}
V_P(h)=S_P(h)+\epsilon \int_0^tf(x)dx -\frac{\epsilon }{2}\frac{f(t)^2}{f'(t)},
   \end{aligned}
  \end{equation}
     where $\epsilon =1$ if $t>0$ and $\epsilon =-1$ otherwise.

 By the assumption $S_P(h)=\frac{2}{3}V_P(h)$, we get from (3.7)
  \begin{equation}\tag{3.8}
  \begin{aligned}
2V_P(h)=6\epsilon \int_0^tf(x)dx -3\epsilon \frac{f(t)^2}{f'(t)}.
   \end{aligned}
  \end{equation}

 After substituting $V_P(h)$ in (3.5) into (3.8), let us differentiate (3.8) with respect to $t$.
  Then we get (3.6). This completes the proof.
 $\square$
\vskip 0.3cm

 Now, it follows from Lemma 3.1 that
   \begin{equation}\tag{3.9}
2\frac{f''(t)}{f'(t)^2}=\frac{tf'(t)-f(t)}{f(t)^2}.
  \end{equation}
Hence we obtain
 \begin{equation}\tag{3.10}
  \begin{aligned}
2(\frac{1}{f'(t)})'=(\frac{t}{f(t)})',
   \end{aligned}
  \end{equation}
which implies that
   \begin{equation}\tag{3.11}
  \begin{aligned}
\frac{2}{f'(t)}=\frac{t}{f(t)}+a,
   \end{aligned}
  \end{equation}
where $a$ is a constant.

After replacing $t$ by $x$,  for $y=f(x)$ we get a differential equation:
\begin{equation}\tag{3.12}
  \begin{aligned}
2ydx-(x+ay)dy=0.
   \end{aligned}
  \end{equation}
  Hence, using a standard method of differential equations, we see that
for some positive constant $b$, $y=f(x)$ satisfies
\begin{equation}\tag{3.13}
  \begin{aligned}
(x-ay)^2=2by.
   \end{aligned}
  \end{equation}
Thus we have
\begin{equation}\tag{3.14}
 \begin{aligned}
 f(x)= \begin{cases}
  \frac{1}{a^2}\{ax+b-\sqrt{2abx+b^2}\}, & \text{if $a\ne0,$} \\
 \frac{x^2}{2b}, & \text{if $a= 0$}.\\
   \end{cases}
  \end{aligned}
  \end{equation}

   Note that
  \begin{equation}\tag{3.15}
  \begin{aligned}
f(0)=f'(0)=0 , f''(0)=\frac{1}{b} \quad \text{and} \quad  f'''(0)=-\frac{3a}{b^2} \quad \text{or} \quad 0.
   \end{aligned}
  \end{equation}
  It follows from (3.13) that  the curve $X$ around an arbitrary point $A_1$
  is an open part of  the parabola defined by
\begin{equation}\tag{3.16}
  \begin{aligned}
x^2-2axy+a^2y^2-2by=0.
   \end{aligned}
  \end{equation}

 \vskip 0.3cm

 Finally using (3.15), in the same manner as in \cite{KK4}, we can show that the curve $X$ is globally an open
 part of a  parabola. This completes the proof of Theorem 1.2.

 \vskip 0.50cm

 \section{Proof of Theorem 1.3}
 \vskip 0.50cm

In this section, we use the main result of \cite{KK4} (Theorem 3 in \cite{KK4}) and Lemma 2.3 in Section 2  in order to prove Theorem 1.3.

 It is obvious that any open part of parabolas satisfy 1) and 2) in Theorem 1.3.

Conversely, suppose that  $X$ denotes a strictly convex $C^{(3)}$ curve  in the plane
 ${\mathbb R}^{2}$ which satisfies  $U_P(h)=\lambda(P)T_P(h)$ for all $P\in X$ and sufficiently small $h>0$.
 Then, it follows from  Lemma 2.1 that $\lambda(P)=\frac{1}{2}$.

First, we note the following which can be easily shown.
\vskip 0.3cm
\noindent {\bf Lemma 4.1.} For a point $P\in X$ and a sufficiently small $h>0$, the following are equivalent.

\noindent 1) $U_P(h)=\frac{1}{2}T_P(h)$,

\noindent 2) $T_P(h)=\frac{1}{2}V_P(h)$,

\noindent 3) $\ell_P(h)=\frac{1}{2}L_P(h)$.
\vskip 0.3cm
Next, using Lemma 2.3 we get the following.
\vskip 0.3cm
\noindent {\bf Lemma 4.2.} Suppose that  $X$ denotes a strictly convex $C^{(3)}$ curve  in the plane
 ${\mathbb R}^{2}$ which satisfies  $U_P(h)=\frac{1}{2}T_P(h)$ for all $P\in X$ and sufficiently small $h>0$.
Then for all $P\in X$ and sufficiently small $h>0$ we have

 \begin{equation}\tag{4.1}
 L_P(h)=\frac{2\sqrt{2}}{\sqrt{\kappa(P)}}\sqrt{h}.
  \end{equation}
\vskip 0.3cm

 \noindent {\bf Proof.} Together with 3) of Lemma 4.1, Lemma 2.3 shows
 \begin{equation}\tag{4.2}
 h\frac{d}{dh}L_P(h)=\frac{1}{2}L_P(h),
  \end{equation}
 which yields for some constant $C=C(P)$
  \begin{equation}\tag{4.3}
L_P(h)=C\sqrt{h}.
  \end{equation}
Thus, Lemma 2.1 completes the proof. $\square$
 \vskip 0.3cm

\vskip 0.3cm

Finally, we prove Theorem 1.3 as follows.

Since $S_P'(h)=L_P(h)$ (\cite{KK4}) and $S_P(0)=0$,
by integrating we get from (4.1)
 \begin{equation}\tag{4.4}
 S_P(h)=\frac{4\sqrt{2}}{3\sqrt{\kappa(P)}}h\sqrt{h}.
  \end{equation}
   Noting $2T_P(h)=hL_P(h)$, one gets from (4.1) and (4.4) that
\begin{equation}\tag{4.5}
 S_P(h)=\frac{4}{3}T_P(h).
  \end{equation}

Theorem 3 of \cite{KK4} states that (4.5) implies $X$ is an open part of a parabola,
completing the proof of Theorem 1.3.
\vskip 0.50cm
 \section{Proof of Theorem 1.4}
 \vskip 0.50cm

In this section, in order to prove Theorem 1.4 we use the main result of \cite{KK4} (Theorem 3 in \cite{KK4}) and Lemma 2.3 in Section 2.

 It is trivial to show that any open part of parabolas satisfy 1) and 2) in Theorem 1.4.

Conversely, suppose that  $X$ denotes a strictly convex $C^{(3)}$ curve  in the plane
 ${\mathbb R}^{2}$ which satisfies $S_P(h)=\lambda(P)W_P(h)$ for all $P\in X$ and sufficiently small $h>0$.
 Then, it follows from  Lemmas 2.1 and 2.2 that $\lambda(P)=\frac{8}{9}$.

First, using Lemma 2.3 we get the following.
\vskip 0.3cm
\noindent {\bf Lemma 5.1.} Suppose that  $X$ denotes a strictly convex $C^{(3)}$ curve  in the plane
 ${\mathbb R}^{2}$ which satisfies  $S_P(h)=\frac{8}{9}W_P(h)$ for all $P\in X$ and sufficiently small $h>0$.
Then for all $P\in X$ and sufficiently small $h>0$ we have

 \begin{equation}\tag{5.1}
 L_P(h)=\frac{2\sqrt{2}}{\sqrt{\kappa(P)}}\sqrt{h}.
  \end{equation}
\vskip 0.3cm

 \noindent {\bf Proof.} Note that
  \begin{equation}\tag{5.2}
 W_P(h)=\frac{1}{2}\{L_P(h)+\ell_P(h)\}h.
  \end{equation}
 Hence, together with Lemma 2.3 the assumption shows that
  \begin{equation}\tag{5.3}
 9S_P(h)=8L_P(h)h-4L_P'(h)h^2.
  \end{equation}

By differentiating (5.3) with respect to $h$ and using $S_P'(h)=L_P(h)$, we get
 \begin{equation}\tag{5.4}
 4L_P''(h)h^2+L_P(h)=0,
  \end{equation}
which is a second order Euler equation.   Its general solutions are given by
\begin{equation}\tag{5.5}
 L_P(h)=C_1\sqrt{h}+C_2\sqrt{h}\ln h,
  \end{equation}
  where $C_1=C_1(P)$ and $C_2=C_2(P)$ are constant.

It follows from Lemma 2.1 that
\begin{equation}\tag{5.6}
 C_1(P)=\frac{2\sqrt{2}}{\sqrt{\kappa(P)}} \quad \text{and } \quad C_2(P)=0.
  \end{equation}
This completes the proof. $\square$
\vskip 0.3cm

Finally, the argument following Lemma 4.2 completes  the proof of Theorem 1.4.
 \vskip 0.5cm

\section{Corollaries}
 \vskip 0.5cm

In this section, we give some corollaries.

Suppose that $X$ is a strictly convex $C^{(3)}$ curve  in the plane
 ${\mathbb R}^{2}$ which  satisfies for all $P\in X$ and sufficiently small $h>0$
  \begin{equation}\tag{6.1}
   \begin{aligned}
 S_P(h)=\lambda(P)V_P(h)^{\mu(P)}.
         \end{aligned}
   \end{equation}
   where $\lambda(P)$ and $\mu(P)$ are some functions.
 Using Lemmas 2.1 and 2.2, by letting $h\rightarrow 0$ we see that
 \begin{equation}\tag{6.2}
   \begin{aligned}
 \lim_{h\rightarrow 0}V_P(h)^{\mu(P)-1}=\frac{2}{3\lambda(P)}.
         \end{aligned}
   \end{equation}
 Since $V_P(h)$ tends to zero as $h\rightarrow 0$, (6.2) shows that
   $\mu(P)=1$.  Hence  we also  obtain from (6.1) that  $\lambda(P)=2/3$.
Thus, from Theorem 1.2 we get
 \vskip 0.3cm

\noindent {\bf Corollary 6.1.} Suppose that  $X$ denotes   a strictly convex  $C^{(3)}$ curve  in the plane ${\mathbb R}^{2}$.
Then, the following are equivalent.
\vskip 0.3cm
\noindent 1) For all $P\in X$ and sufficiently small $h>0$,
 $X$ satisfies $ S_P(h)=\lambda(P)V_P(h)^{\mu(P)}$, where  $\lambda(P)$ and $\mu(P)$ are some functions.

\noindent 2) For all $P\in X$ and sufficiently small $h>0$, $X$ satisfies $ S_P(h)=\frac{2}{3}V_P(h)$.

\noindent 3) $X$ is an open part of a parabola.

\vskip 0.3cm
The similar argument as in the proof of Corollary 6.1 shows the following.
\vskip 0.3cm

\noindent {\bf Corollary 6.2.} Suppose that  $X$ denotes   a strictly convex $C^{(3)}$ curve  in the plane ${\mathbb R}^{2}$.
Then, the following are equivalent.
\vskip 0.3cm
\noindent 1) For all $P\in X$ and sufficiently small $h>0$,
 $X$ satisfies $ U_P(h)=\lambda(P)T_P(h)^{\mu(P)}$, where  $\lambda(P)$ and $\mu(P)$ are some functions.

\noindent 2) For all $P\in X$ and sufficiently small $h>0$, $X$ satisfies $ U_P(h)=\frac{1}{2}T_P(h)$.

\noindent 3) $X$ is an open part of a parabola.

\vskip 0.3cm

\noindent {\bf Corollary 6.3.} Suppose that  $X$ denotes   a strictly convex $C^{(3)}$ curve  in the plane ${\mathbb R}^{2}$.
Then, the following are equivalent.
\vskip 0.3cm
\noindent 1) For all $P\in X$ and sufficiently small $h>0$,
 $X$ satisfies $ S_P(h)=\lambda(P)W_P(h)^{\mu(P)}$, where  $\lambda(P)$ and $\mu(P)$ are some functions.

\noindent 2) For all $P\in X$ and sufficiently small $h>0$, $X$ satisfies $ S_P(h)=\frac{8}{9}W_P(h)$.

\noindent 3) $X$ is an open part of a parabola.

\vskip 0.3cm

\vskip 1.0 cm

Department of Mathematics, \par Chonnam National University,\par
Kwangju 500-757, Korea

{\tt E-mail:dosokim@chonnam.ac.kr} \vskip 0.3 cm

Department of Mathematics, \par Chonnam National University,\par
Kwangju 500-757, Korea

{\tt E-mail:yong4625@naver.com} \vskip 0.3 cm

Department of Mathematics, \par Kyungpook National University,\par
Taegu 702-701, Korea

{\tt E-mail:yhkim@knu.ac.kr} \vskip 0.3 cm

Department of Mathematics, \par Chonnam National University,\par
Kwangju 500-757, Korea

{\tt E-mail:dhpark3331@chonnam.ac.kr} \vskip 0.3 cm

\end{document}